\pdfoutput=1
\documentclass{article}
\usepackage[margin=25mm]{geometry}
\usepackage{amsmath}
\usepackage{amsfonts}
\usepackage{amssymb}
\usepackage{graphicx}

\usepackage[utf8]{inputenc} 
\usepackage[T2A]{fontenc}
\usepackage[english]{babel}
\usepackage[unicode, pdftex]{hyperref}

\usepackage{xcolor}
\usepackage[all]{xy}

\newtheorem{theorem}{Theorem}

\title{Optimal Reeb graphs on  two- and three-connected planar polygon}
\author{Oleksandr Pryshliak, Karolina Haieva}

\date{\today}

\begin{document}
\maketitle

\begin{abstract}
To investigate the topological structure of planar polygon decomposition on trapezoids, which is formed by height functions. We use the oriented  Reeb graph of the function with a marked vertex. We describe all possible optimal Reeb graphs in the case of polygon in general position with one local minimum and one local maximum.  By optimal Reeb graph we mean Reeb graph, which cann't be obtaned from other Reeb graph by subdivision of an edge or a leaf attaching. In this case polygon is a triangle with triangle holes. Constructed Reeb graphs give topological structures of trapezoid maps on two-connected polygon with six vertexes and three-connected polygon with nine vertexes.
\end{abstract} \hspace{10pt}



\section*{Introduction}

In 1946, Reeb introduced \cite{Reeb1946} a graph, which describes the topological properties of functions on two-dimensional manifolds. In the case of simple Morse functions on closed oriented two-dimensional manifolds, the  Reeb graph is a complete topological invariant of fiberwise equivalence. It becomes a complete topological invariant under topological equivalence, if a linear order is specified on the set of vertices. In the case of a non-orientable two-dimensional manifold \cite{lychak2009morse}, as well as manifolds with boundary \cite{hladysh2017topology, hladysh2019simple}, additional information is needed to construct a complete topological invariant. Reeb graph on non-compact 2-manifolds wos investigated in \cite{Kronrod1950, prishlyak2002morse}. In this case it can be a non-Hausdorff space. If the function has by  several critical points on its level, then we need the information about the structure of the function in the neighborhood of critical levels \cite{Bolsinov2004}.

Another way to specify the topological structure of Morse functions is to use Morse-Smale vector fields with given function values at singular points \cite{lychak2009morse, Smale1961}. Therefore, the topological classification of Morse--Smale vector fields is closely related to the classification of functions.

Graphs as topological invariant of functions was used in \cite{prishlyak2001conjugacy, hladysh2019simple, hladysh2017topology,  prishlyak2002morse, Prishlyak2000,   lychak2009morse, prish2002Morse, prish2015top, prish1998sopr,  Bilun2002,  Sharko1993}, on surfaces with the boundary in \cite{Hladysh2016, hladysh2019simple} and on closed 3-manifolds in  \cite{prishlyak1999equivalence, kkp2013, prishlyak2003regular}. 

By analogy to Reeb graph, complete invariants of  flows ar also graph (destiguioshed graphs). Its was constructed in \cite{Kybalko2018,  prishlyak1997graphs,   Prishlyak2021,   Kybalko2018} on classifications of flows on closed 2- manifolds and 
\cite{prishlyak2003sum} on 2-manifolds with the boundary and in
\cite{prish1998vek,  prish2001top, prishlyak2003regular, Prishlyak2002b,  Hatamian2020} as Heegard diagrams of 3-manifolds.

We recoment \cite{prishlyak1997graphs, Harary69, pontr86, tatt88, HW68, GT87}, where you can find the main invariants of graphs and their embeddings in surfaces.

Effective computing of Reeb graphs was done in \cite{Doraiswamy2013}.

The main purpose of this paper is to find all possible Reeb graphs of height function with minimal number of vertexes on two- and three-connected polygons.  


\textbf{The main results.}

In Chapter 1, we describe the relationship between Reeb graphs and trapezoidal maps of polygons.

Chapter 2 describes the properties of Reeb graphs on simple polygons and lists all graphs with no more than 5 vertices.

Chapter 3 describes the properties of the optimal Reeb graph. In particular, theorem 1 proves that such a graph has 3n vertices on an n-connected polygon, and theorem 2 describes the number of vertices of each degree.

Chapter 4 lists all optimal graphs on biconnected polygons.

Chapter 5 describes all optimal graphs on three-connected polygons
 
\section{Trapezoidal map and Reeb graph}
 



Let $L$ be a polygon with a set of vertices $V=\{v_i(x_i,y_i), 1\le i\le m\}$ and edges $E=\{e_j=\{v_j,w_j\}, v_j \in V,w_j\in V\}, 1\le j\le n$. The polygon is \textit{in general position} with the height function $h(x,y)=y$, if $v_i\ne v_j$ implies $y_i \ne y_j$. In the future, we will consider polygons that are in general position with the function of height.
 
By \textit{trapezoid }we will understand one of the following polygons: trapeze, parallelogram, triangle.

Let's choose a sufficiently large rectangle $S=[a,b] \times[c,d]$, that $V \subset \text{Int}\ S$, (for example, $a=\min x_i - 1, \ b= \max x_i + 1, \ c=\min y_i - 1,\ d=\max y_i + 1$).

For each vertex $v_i$, we draw a horizontal ray to the left to the first intersection with another edge of $L$ or the line $x=a$. Similarly, we draw the ray to the right to the first point of intersection with the edges of $L$ or the line $x=b$. Consider those segments that lie inside the polygon $L$. They split $L$ into trapezoids. The partition of $S$ into trapezoids constructed in this way is called a trapezoidal map. We will denote it by $T(L)$. By construction, a trapezoidal map is a flat rectilinear graph. Since no more than 2 horizontal segments emerge from each vertex of the polygon $L$, and none from the lowest and highest vertices, the total number of vertices of the trapezoidal map does not exceed $3n-4$.

\begin{figure}[ht]
\center{ \includegraphics[width=0.4\linewidth]{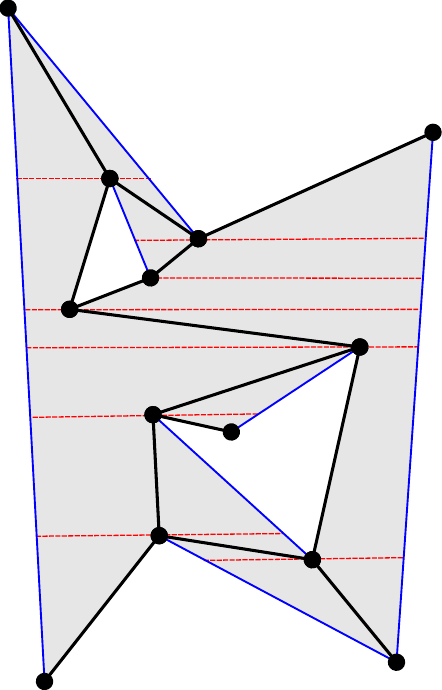} }
\caption{Example of the trapezoidal map (with red dashed bases) and Reeb graph (black)}
\label{Trap}
\end{figure}

A component of the level line $h^{ -1} (y)$ of the height function is called \textit{a fiber}. 

Two  functions are called \textit{fibre equivalent} if there is a homeomorphism of the polygon onto the polygon which maps fibers of one function to fibers of another one. 

The quotient space $M/ \sim$ with orientation of edges according to the direction of the increase of the function is called a \textit{Reeb graph}, where $h : M \to R$ is the height function, $x_1 \sim x_2$ if $x_1$ and $x_2$ belongs to one fiber. Reeb graphs are considered to within the isomorphism of oriented graphs. 
Examples of Reeb graphs are shown on Fig. \ref{Trap}.

By construction we have one to one correspondence trapezoids to edges
of Reeb graph. Thus, the Reeb graph define the structure of the trapezoidal map. 

In addition, vertexes of degree 2 are divided on left and right. Vertex is \textit{left} if right horizontal ray intersect the polygon in its small neighbourhood. Otherwise vertex is \textit{right}.

Two simple height functions on the orientable surface are \textit{fiber equivalent} if and only if their Reeb graphs are isomorphic ([1], Theorem 2.4, p. 71). On the non-orientable surfaces there are firer non-equivalent functions with the same Reeb graphs. So additional invariants are nessersity to clasify height functions. 
But in the case of planar polygon in general position Reeb graph is complete invariant of height functions and structure of the trapezoidal map.

\section{Reeb graph on simple polygon}

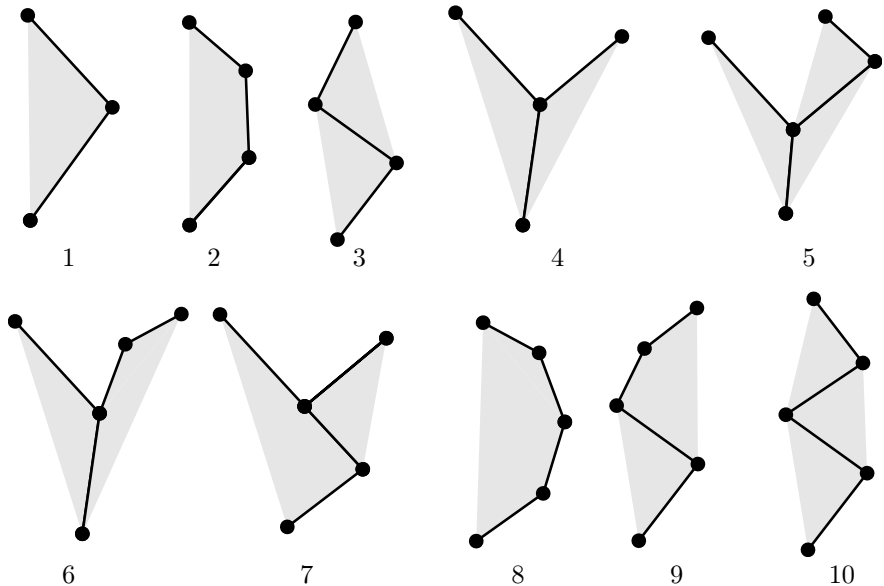
\begin{figure}[ht]
\center{ \includegraphics[width=0.7\linewidth]{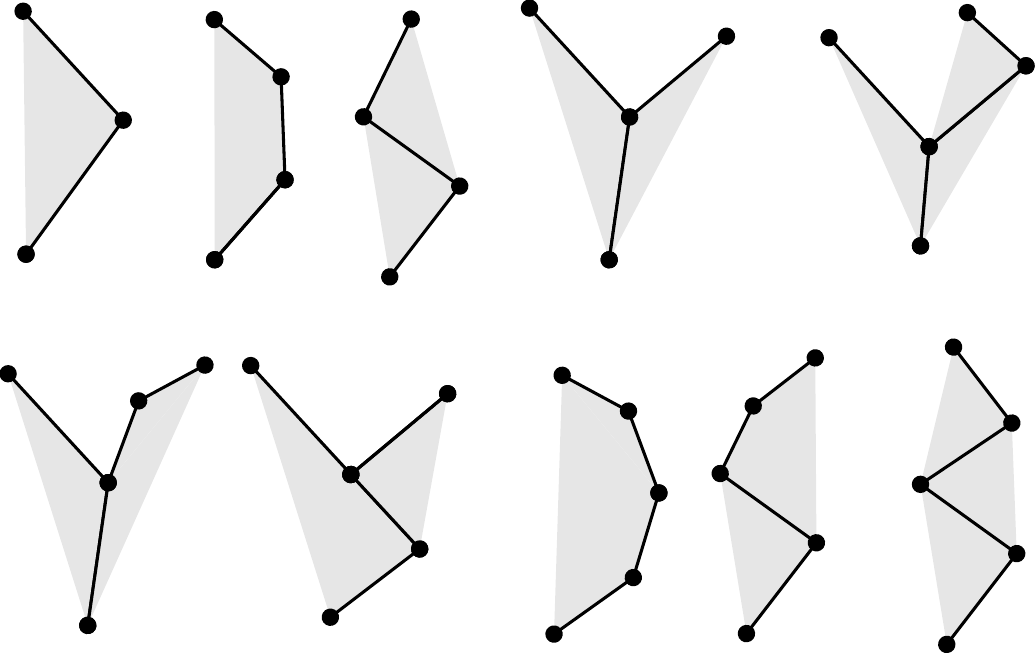} }
\put (-310,110) {1}
\put (-255,110) {2}
\put (-200,110) {3}
\put (-125,110) {4}
\put (-30,110) {5}
\put (-310,-10) {6}
\put (-220,-10) {7}
\put (-140,-10) {8}
\put (-80,-10) {9}
\put (-20,-10) {10}
\caption{Reeb graphs   on simple polygon with 3, 4 and 5 vertexes}
\label{Morse10}
\end{figure}

By the definition of the Reeb graph, it is homotopically equivalent to a polygon. So, if a polygon is simple, which means simply-connected, then the Reeb graph has the same properties. This means that it is a tree. Since the simple polygon with the smallest number of vertices is a triangle, the Reeb graph has at least three vertices. This graph is shown in fig.
\ref{Morse10}.1.

Trees with 4 vertexes are either chains, two possible options are shown in Fig. \ref{Morse10}.2,3, or the tree is like the letter Y - \ref{Morse10}.4.
A tree with 5 vertices is separated from a tree with 4 vertices by dividing one of the edges. All possible options to divide Y are shown in Fig. \ref{Morse10}.5--7. Chains  with 5 vertexes are shown in Fig. \ref{Morse10}.8-10.

\section{Optimal Reeb graph}







The Reeb graph $G$ of a polygon is called optimal if there is no other Reeb graph $H$ of a polygon such that $G$ and $H$ are homotopically equivalent and the number of vertices of the graph $H$ is less than the number of vertices of the graph $G$.

This means that the Reeb graph is not optimal, if  it is possible to subtract from another Reeb graph by the edge subdivision or add a leaf.

The polygon is called \textit{n-connected} if it is connected, and  its boundary has an n component of connectivity.

\begin{theorem} The optimal Reeb graph on an n-connected polygon has 3n vertices.

\end{theorem}
\textbf{Proof.} Each component of the boundary is a polygon which has at least 3 vertexes. Therefore, the number of vertices is no less than 3n. Let us show that there is a height function that its  Reeb graph has 3n vertices. For this purpose, let us look at the rather large triangle A, the vertexes of which lie at different heights. Consider $n-1$ triangle $A_i \subset \text{Int} A$, which does not have intersection. The polygon $B=A \setminus \cup_{i=1}^{n-1}A_i$ contains 3n vertices, and the number of vertices of Reeb graph is same.

The completed theorem shows a method for generating all optimal Reeb graphs.

\begin{theorem}
The optimal Reeb graph of the height function on the n-connected polygon has the following properties:

1) the graph is  n-connected planar graph;

2) the graph has n vertexes of degree 2, two vertexes of degree 1, the other vertices have degree 1 or 3;

3) only vertices of degree 1 are sinks and sources.
\end{theorem}
\textbf{Proof.} 1) It followed from the definition of Reeb graph.

2) Each triangle has a regular vertex of degree 2. The height function has a maximum and a minimum. If there is other vertex of degree one or two, then the graph is not optimal.

3) sinks and sources correspond to local extrema of the function, and vertices of degree 2 and 3 correspond to saddle points.
\section{Optimal Reeb graph on two-connected polygon}
As it was meant above all, Reeb graph of two-connected polygon is optimal, if it has 6 vertices, and the polygon has a triangle with  a hole that looks like a triangle. 
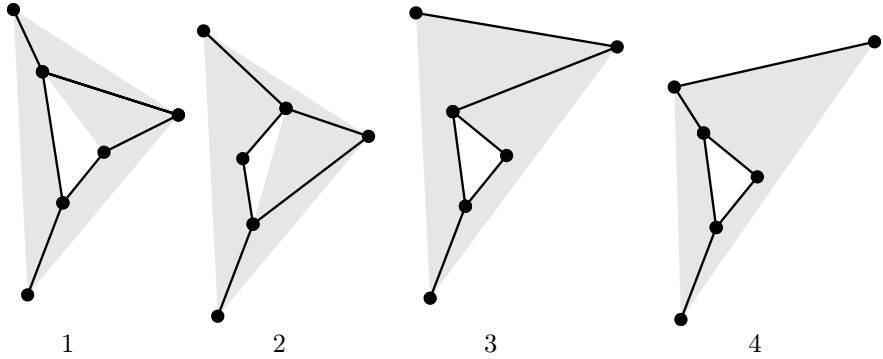
\begin{figure}[ht]
\center{ \includegraphics[width=0.7\linewidth]{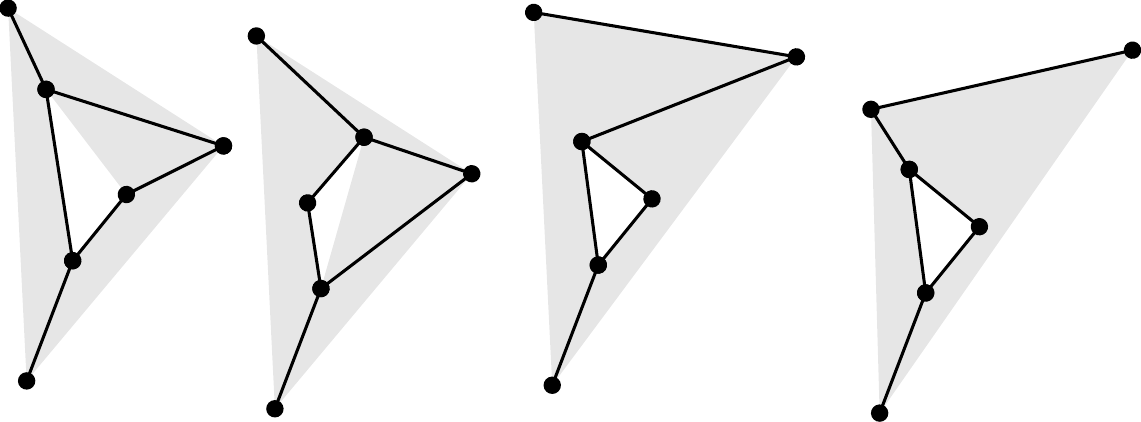} }
\put (-310,-10) {1}
\put (-230,-10) {2}
\put (-150,-10) {3}
\put (-50,-10) {4}
\caption{Optimal Reeb graph of two-connected polygon}
\label{Morse12}
\end{figure}

The \textit{middle vertex} of the triangle is called the vertex that is not the lowest or the highest for the triangle. Since the axial symmetry of the vertical axis is identical, we can take into account without limiting the formality that the middle vertex  on the external triangle is the rightmost (which is the largest abscise). So, the axial symmetry of exactly horizontal straight lines leads to equivalent graphs. Therefore, we can take into account that at least two vertexes of the inner triangle lie lower than the middle vertex of the outer triangle. There is two or three such vertexes. In addition, the inner triangle can be left (middle apex of the left) or right (right). Therefore, there are 4 different possibilities for the development of the middle triangle and according to the different Reeb graphs. Its is shown in Fig. 3.

Summarizing the results obtained above, we have the following 
\begin{theorem} 
There are 4 non-isomorphic optimal Reeb graphs on two-connected  polygons.
\end{theorem}

\section{Optimal Reeb graph of three-connected polygon}
We denote by ABC the outer triangle, and by DEF and GHI  the inner triangles, such that A>B>C, D>E>F, G>H>I, D>G.

From now on, as before, we will be able to see only such triangles, in which B is the right vertex of the outer triangle and no less than three vertices of the inner triangles lie below it.
We call the inner triangle DEF by upper and  we  call GHI by the lower triangle. Let's look at it first if F>G.
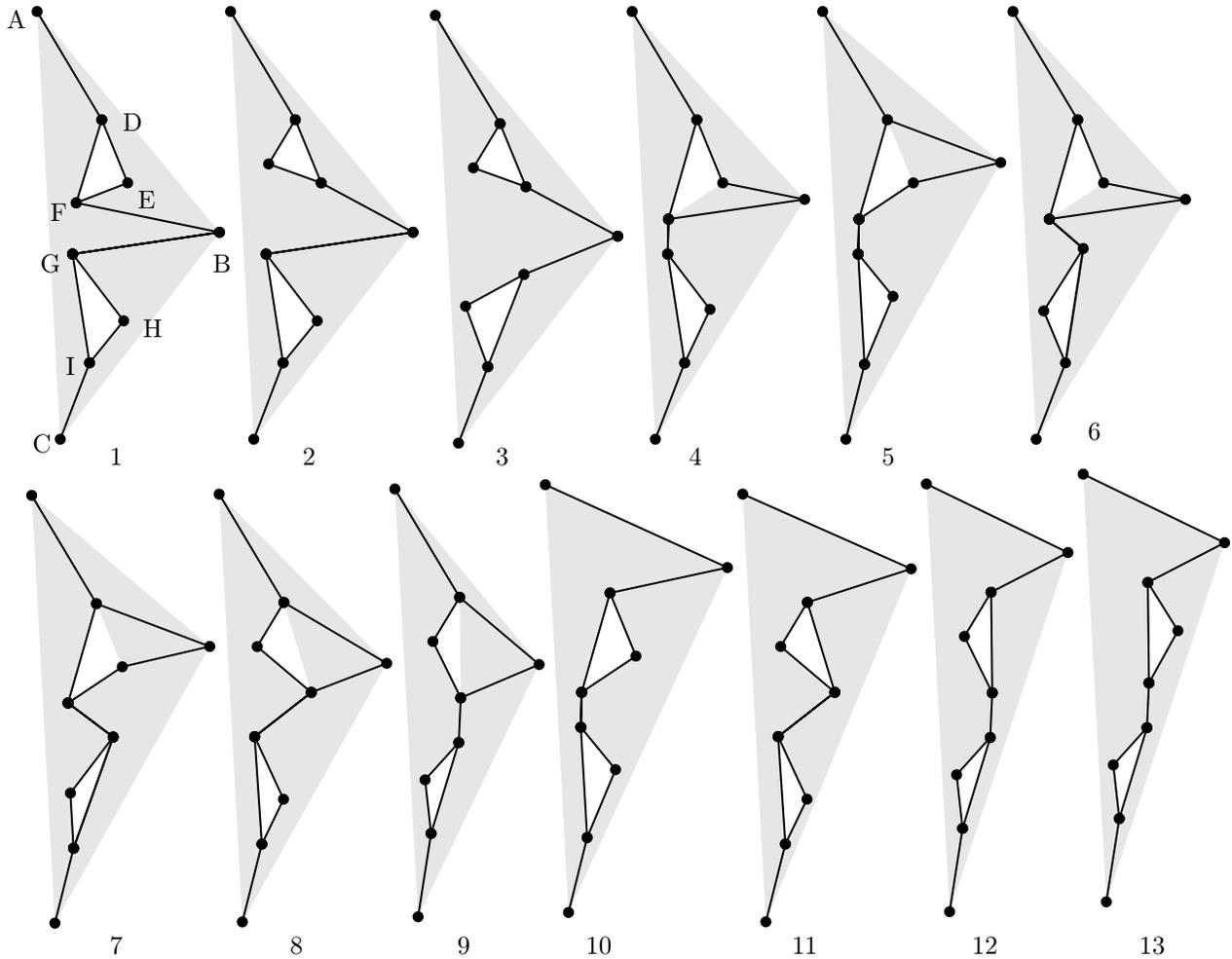
\begin{figure}[ht]
\center{ \includegraphics[width=0.99\linewidth]{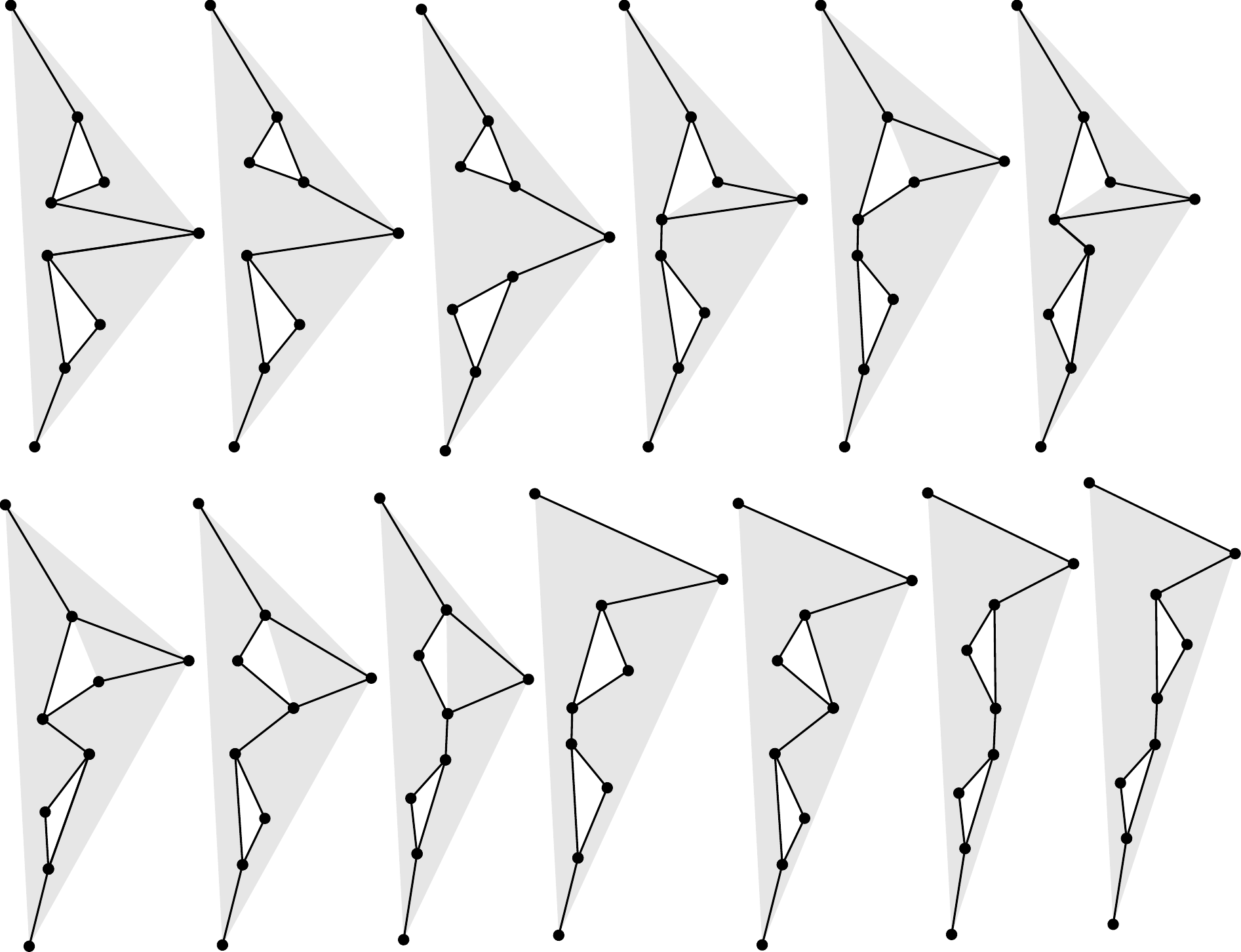} }
\put (-475,350) {A}
\put (-395,255) {B}
\put (-465,185) {C}
\put (-430,310) {D}
\put (-424,280) {E}
\put (-458,275) {F}
\put (-462,255) {G}
\put (-422,230) {H}
\put (-452,215) {I}
\put (-435,180) {1}
\put (-360,180) {2}
\put (-285,180) {3}
\put (-210,180) {4}
\put (-135,180) {5}
\put (-55,190) {6}
\put (-435,-10) {7}
\put (-365,-10) {8}
\put (-300,-10) {9}
\put (-250,-10) {10}
\put (-170,-10) {11}
\put (-100,-10) {12}
\put (-35,-10) {13}
\caption{Optimal Reeb graph of three-connected polygon. Part 1}
\label{Morse3}
\end{figure}

We suppose F>B>G.  There are three options: 1) the inner triangles are  right Fig. \ref{Morse3}.1, 2) one of them is left, and the other is right Fig. \ref{Morse3}.2, 3) the inner triangle are left Fig. \ref{Morse3}.3. By moving the middle outer vertex up, a new expansion of the triangle is removed. If we use it for Fig.\ref{Morse3}-1, we obtain Fig. \ref{Morse3}.4, Fig. \ref{Morse3}.5 and Fig. \ref{Morse3}.10. If we use it for Fig. \ref{Morse3}.2 we obtain Fig. \ref{Morse3},8 and \ref{Morse3}.11. If we use it for horizontal symmetric of Fig. \ref{Morse3}.2, we obtain Fig. \ref{Morse3}6, 7,13. 
Figure \ref{Morse3}.3 gives Fig. \ref{Morse3}.9,12.



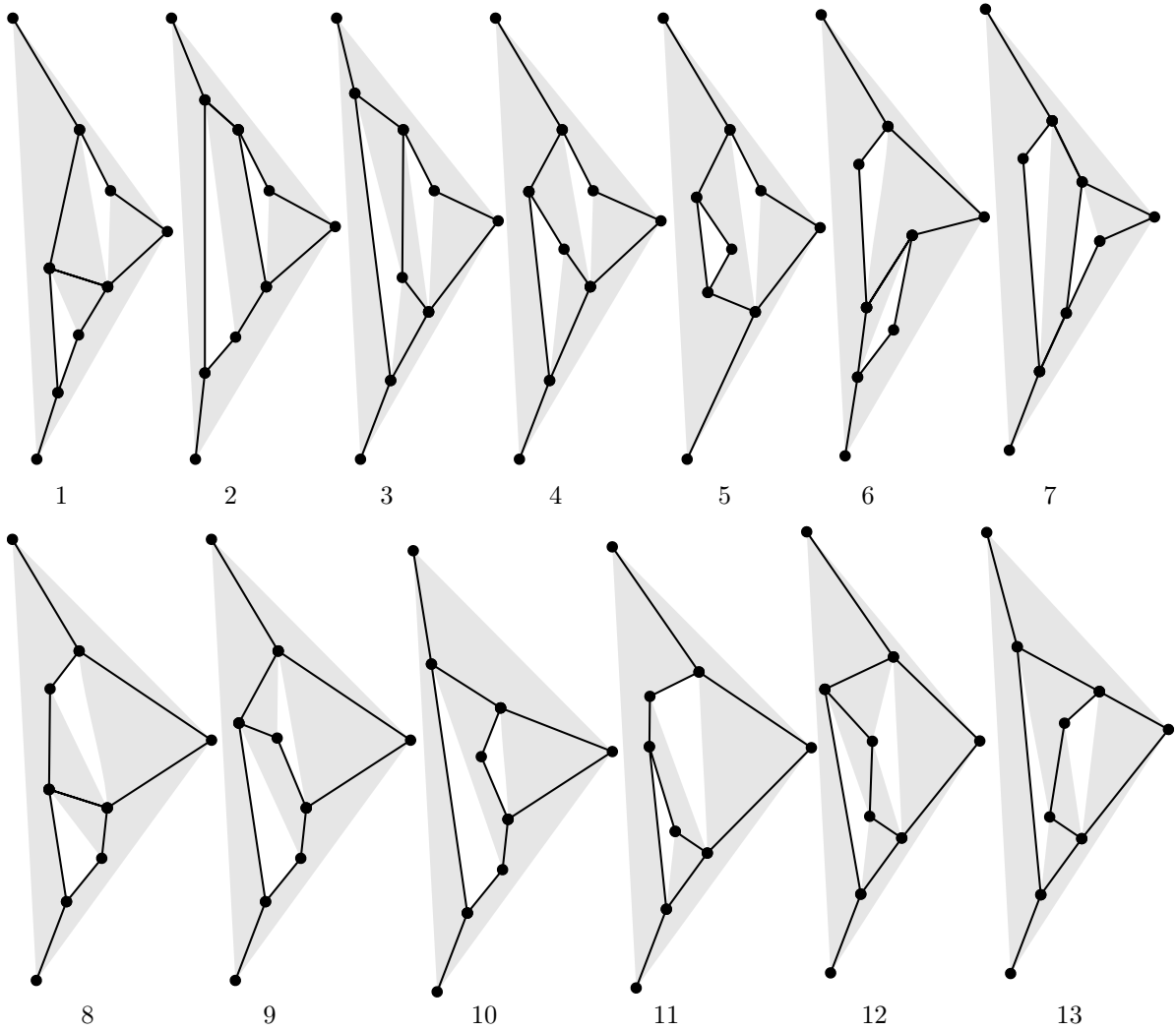
\begin{figure}[ht]
\center{ \includegraphics[width=0.95\linewidth]{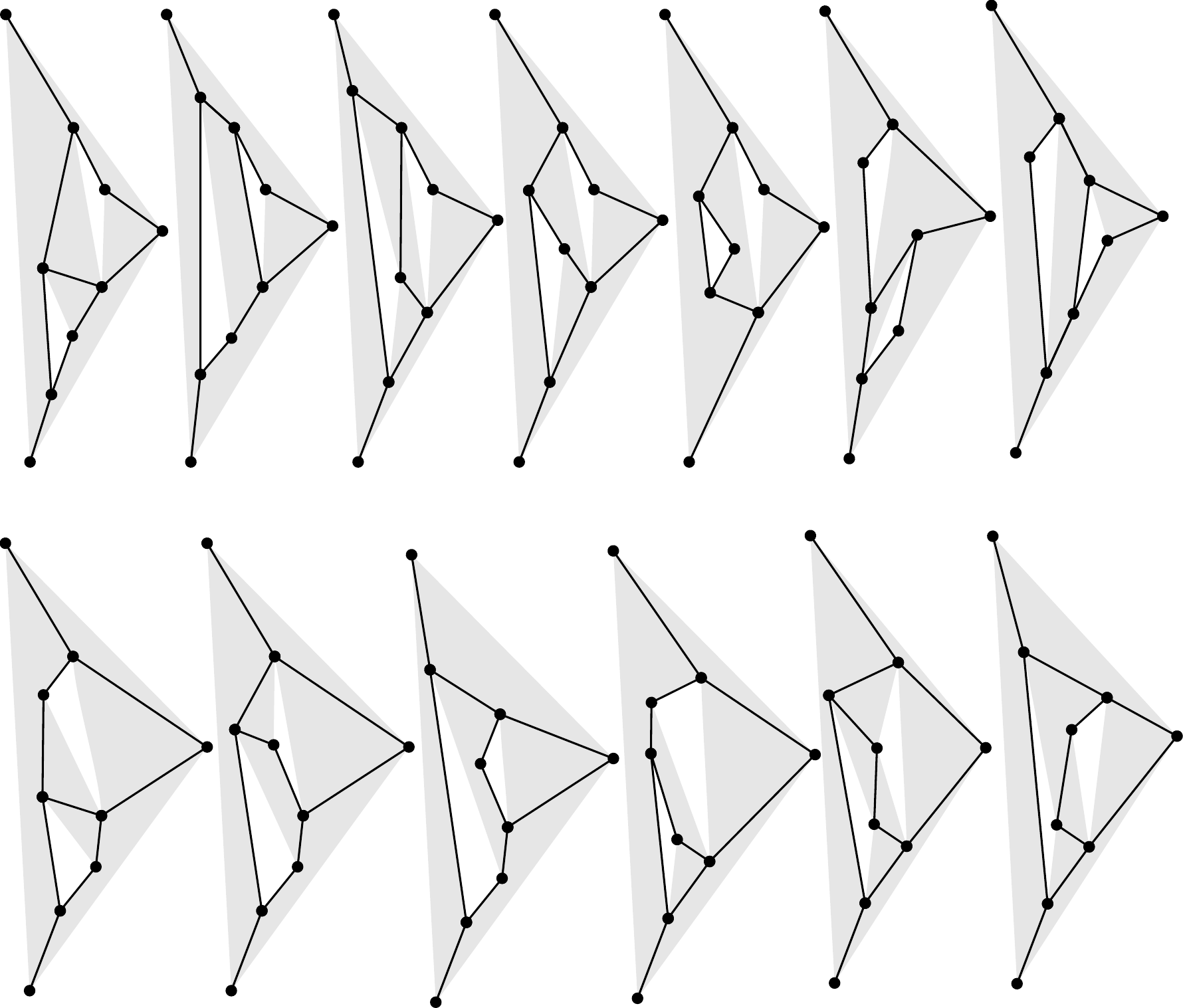} }
\put (-430,190) {1}
\put (-365,190) {2}
\put (-305,190) {3}
\put (-240,190) {4}
\put (-175,190) {5}
\put (-120,190) {6}
\put (-50,190) {7}
\put (-420,-10) {8}
\put (-350,-10) {9}
\put (-270,-10) {10}
\put (-200,-10) {11}
\put (-120,-10) {12}
\put (-45,-10) {13}
\caption{Optimal Reeb graph of three-connected polygon. Part 2
}
\label{Morse4}
\end{figure}




Next, we consider the cases when D>G>F. Fig. {Morse4} shows all possible arrangements of internal triangles. In the case of two solutions of a pair of internal triangles, which are symmetrical about the vertical axis, we depict only one of them. Therefore, the vertex B can be located both on the right and on the left.

In Fig. \ref{Morse4}.1 we have D>E>G>F>H>I. Possible positions of the vertex B on the left - 4, and on the right - 6.

In Fig. \ref{Morse4}.2 we have D>G>H>I>E>F. Possible positions of the vertex B on the left - 3, and on the right - 7.

In Fig. \ref{Morse4}.3 we have D>G>H>E>I>F. Possible positions of the vertex B on the left - 3, and on the right - 6.

In Fig. \ref{Morse4}.4, possible positions of the vertex B on the left - 4, and on the right - 5.

In Fig. \ref{Morse4}.5, possible positions of the vertex B on the left - 5, and on the right - 4.

In Fig. \ref{Morse4}.6, possible positions of the vertex B on the left - 5, and on the right - 5.

In Fig. \ref{Morse4}.7, possible positions of the vertex B on the left - 4, and on the right - 6.

In Fig. \ref{Morse4}.8, possible positions of the vertex B on the left - 5, and on the right - 5.

In Fig. \ref{Morse4}.9, possible positions of the vertex B on the left - 4, and on the right - 5.

In Fig. \ref{Morse4}.10, possible positions of the vertex B on the left - 3, and on the right - 6.

In Fig. \ref{Morse4}.11, possible positions of the vertex B on the left - 5, and on the right - 4.

In Fig. \ref{Morse4}.12, possible positions of the vertex B on the left - 4, and on the right - 4.

In Fig. \ref{Morse4}.13, possible positions of the vertex B on the left - 3, and on the right - 5.

Summarizing the results obtained above, we have the following 
\begin{theorem} 
There are 133 non-isomorphic optimal Reeb graphs on three-connected  polygons.
\end{theorem}

\section*{Conclusion} 


In this paper, we have found all possible structures of optimal Reeb graphs on two-connected (4 structures) and three-connected (133 structures) polygons. Any other Reeb graph can be obtained from the optimal one by using edge splitting and leaf addition operations. The resulting Reeb graphs describe the structures of trapezoidal maps on these polygons.


\begin{thebibliography}{10}




\bibitem{Bilun2002}
S.~Bilun and A.~Prishlyak.
\newblock The closed morse 1-forms on closed surfaces.
\newblock {\em Visn., Mat. Mekh., Kyv. Univ. Im. Tarasa Shevchenka},
  2002(8):77--81, 2002.



\bibitem{Bolsinov2004}
A.V. Bolsinov and A.T. Fomenko.
\newblock {\em Integrable Hamiltonian systems. Geometry, Topology,
  Classification}.
\newblock A CRC Press Company, Boca Raton London New York Washington, D.C.,
  2004.
\newblock 724 p.


\bibitem{Doraiswamy2013}
H. Doraiswamy, V. Natarajan
\newblock {\em Computing Reeb Graphs as a Union of Contour Trees}.
\newblock IEEE Trans Vis Comput Graph
. 2013 Feb;19(2):249-62. doi: 10.1109/TVCG.2012.115.

\bibitem{GT87}
J.~Gross and T.~Tucker.
\newblock {\em Topological graph theory}.
\newblock Wiley, 1987.

\bibitem{Harary69}
F.~Harary.
\newblock {\em Graph Theory}.
\newblock Westview Press, 1969.

\bibitem{Hatamian2020}
C. Hatamian and A.~Prishlyak.
\newblock Heegaard diagrams and optimal morse flows on non-orientable
  3-manifolds of genus 1 and genus 2.
\newblock {\em Proc. Int. Geom. Cent.}, 13(3):33--48, 2020.
\href{http://dx.doi.org/10.15673/tmgc.v13i3.1779}{\path{doi: 10.15673/tmgc.v13i3.1779}}.

\bibitem{HW68}
P.~Hilton and S.~Wylie.
\newblock {\em Homology theory: An introduction to algebraic topology}.
\newblock Cambridge University Press, 1968.

\bibitem{Hladysh2016}
B.~I. Hladysh and A.~O. Pryshlyak.
\newblock Functions with nondegenerate critical points on the boundary of the
  surface.
\newblock {\em Ukrainian Mathematical Journal}, 68(1):29--40, jun 2016.
\href{http://dx.doi.org/10.1007/s11253-016-1206-5}{\path{doi: 10.1007/s11253-016-1206-5}}.

\bibitem{hladysh2017topology}
B.I. Hladysh and A.O. Prishlyak.
\newblock Topology of functions with isolated critical points on the boundary
  of a 2-dimensional manifold.
\newblock {\em SIGMA. Symmetry, Integrability and Geometry: Methods and
  Applications}, 13:050, 2017.
\href{http://dx.doi.org/0.3842/SIGMA.2017.050}{\path{doi: 0.3842/SIGMA.2017.050}}.

\bibitem{hladysh2019simple}
B.I. Hladysh and A.O. Prishlyak.
\newblock Simple morse functions on an oriented surface with boundary.
\newblock {\em Журнал математической физики,
  анализа, геометрии}, 15(3):354--368, 2019.
\href{http://dx.doi.org/10.15407/mag15.03.354}{\path{doi: 10.15407/mag15.03.354}}.

\bibitem{Kronrod1950}
A.S. Kronrod.
\newblock Functions of two variables.
\newblock {\em Russian Mathematical Surveys}, 5:24--134, 1950.

\bibitem{Kybalko2018}
Z.~Kybalko, A.~Prishlyak, and R.~Shchurko.
\newblock {Trajectory equivalence of optimal Morse flows on closed surfaces}.
\newblock {\em {Proc. Int. Geom. Cent.}}, 11(1):12--26, 2018.
\href{http://dx.doi.org/10.15673/tmgc.v11i1.916 }{\path{doi: 10.15673/tmgc.v11i1.916}}.


\bibitem{lychak2009morse}
D.P. Lychak and A.O. Prishlyak.
\newblock Morse functions and flows on nonorientable surfaces.
\newblock {\em Methods of Functional Analysis and Topology}, 15(03):251--258,
  2009.




\bibitem{Prishlyak2021}
A.~Prishlyak, A.~Prus, and S.~Guraka.
\newblock Flows with collective dynamics on a sphere.
\newblock {\em Proc. Int. Geom. Cent.}, 14(1):61--80, 2021.
\href{http://dx.doi.org/10.15673/tmgc.v14i1.1902}{\path{doi: 10.15673/tmgc.v14i1.1902}}.

\bibitem{prishlyak1997graphs}
A.O. Prishlyak.
\newblock On graphs embedded in a surface.
\newblock {\em Russian Mathematical Surveys}, 52(4):844, 1997.
\href{http://dx.doi.org/10.1070/RM1997v052n04ABEH002074}{\path{doi: 10.1070/RM1997v052n04ABEH002074}}.

\bibitem{prishlyak1999equivalence}
A.O. Prishlyak.
\newblock Equivalence of morse function on 3-manifolds.
\newblock {\em Methods of Func. Ann. and Topology}, 5(3):49--53, 1999.

\bibitem{Prishlyak2000}
A.O. Prishlyak.
\newblock Conjugacy of morse functions on surfaces with values on a straight
  line and circle.
\newblock {\em Ukrainian Mathematical Journal}, 52(10):1623--1627, 2000.
\href{http://dx.doi.org/10.1023/A:1010461319703}{\path{doi: 10.1023/A:1010461319703}}.

\bibitem{prishlyak2001conjugacy}
A.O. Prishlyak.
\newblock Conjugacy of morse functions on 4-manifolds.
\newblock {\em Russian Mathematical Surveys}, 56(1):170--171, 2001.
\href{http://dx.doi.org/10.1070/RM2001v056n01ABEH000370}{\path{doi: 10.1070/RM2001v056n01ABEH000370}}.


\bibitem{prishlyak2002morse}
A.O. Prishlyak.
\newblock Morse functions with finite number of singularities on a plane.
\newblock {\em Meth. Funct. Anal. Topol}, 8:75--78, 2002.



\bibitem{Prishlyak2002b}
A.O. Prishlyak.
\newblock Topological equivalence of smooth functions with isolated critical
  points on a closed surface.
\newblock {\em Topology and its Applications}, 119(3):257--267, 2002.
\href{http://dx.doi.org/10.1016/S0166-8641(01)00077-3}{\path{doi: 10.1016/S0166-8641(01)00077-3}}.

\bibitem{prishlyak2003sum}
A.O. Prishlyak.
\newblock On sum of indices of flow with isolated fixed points on a stratified
  sets.
\newblock {\em Zhurnal Matematicheskoi Fiziki, Analiza, Geometrii [Journal of
  Mathematical Physics, Analysis, Geometry]}, 10(1):106--115, 2003.

\bibitem{prishlyak2003regular}
A.O. Prishlyak.
\newblock Regular functions on closed three-dimensional manifolds.
\newblock {\em Dopov. Nats. Akad. Nauk Ukr. Mat. Prirodozn. Tekh. Nauki},
  8:21--24, 2003.




\bibitem{Reeb1946}
G.~Reeb.
\newblock Sur les points singuliers d'une forme de pfaff complétement
  intégrable ou d'une fonction numérique.
\newblock {\em C.R.A.S. Paris}, 222:847—849, 1946.

\bibitem{Sharko1993}
V.V. Sharko.
\newblock {\em Functions on manifolds. Algebraic and topological aspects.},
  volume 131 of {\em Translations of Mathematical Monographs}.
\newblock American Mathematical Society, Providence, RI, 1993.

\bibitem{Smale1961}
S.~Smale.
\newblock On gradient dynamical systems.
\newblock {\em Ann. of Math.}, 74:199--206, 1961.


\bibitem{kkp2013}
В.М. Кузаконь, В.Ф. Кириченко, and О.О. Пришляк.
\newblock Гладкi многовиди. Геометричнi та
  топологiчнi аспекти.
\newblock {\em Працi Iнcтитуту математики НАН
  України.—2013.—97.—500 с}, 2013.


\bibitem{pontr86}
Л.~С. Понтрягин.
\newblock {\em Основы комбинаторной топологии}.
\newblock Наука, 1986.

\bibitem{prish1998vek}
А.О. Пришляк.
\newblock Векторные поля Морса--Смейла с
  конечным числом особых траекторий на
  трехмерных многообразиях.
\newblock {\em Доповiдi НАН України}, 6:43--47, 1998.

\bibitem{prish1998sopr}
А.О. Пришляк.
\newblock Сопряженность функций Морса.
\newblock {\em Некоторые вопросы совр.
  математики. Институт математики АН
  Украины, Киев}, 1998.

\bibitem{prish2001top}
А.О. Пришляк.
\newblock Топологическая эквивалентность
  функций и векторных полей Морса—Смейла
  на трёхмерных многообразиях.
\newblock {\em Топология и геометрия. Труды
  Украинского мат. конгресса}, pages 29--38, 2001.


\bibitem{prish2002Morse}
О.О. Пришляк.
\newblock {\em Теорія Морса}.
\newblock Київський університет, 2002.

\bibitem{prish2015top}
О.О. Пришляк.
\newblock {\em Топологія многовидів}.
\newblock Київський університет, 2015.

\bibitem{tatt88}
У. Татт.
\newblock {\em Теория графов}.
\newblock Мир, 1988.


\end{thebibliography}

\textsc{Taras Shevchenko National University of Kyiv}

Oleksandr Pryshliak \ \textit{Email address:} \text{ prishlyak@knu.ua} \ \ \ \
\textit{ Orcid ID:} \text{0000-0002-7164-807X}

Karolina Haieva \ \ \textit{Email address:} \text{ gayevacarolina@knu.ua} \ 
\textit{ Orcid ID:} \text{0009-0004-2613-4055}

\end{document}